\documentclass[a4paper,11pt]{amsart}

\usepackage{graphicx}
\usepackage{mathptmx}
\usepackage{amsmath}
\usepackage{amssymb}
\usepackage{enumitem}
\usepackage{xcolor}
\usepackage{pgfplots}

\newmuskip\pFqmuskip

\newcommand*\pFq[6][8]{%
  \begingroup 
  \pFqmuskip=#1mu\relax
  \mathcode`=\string"8000
  \begingroup\lccode`\~=`\,
  \lowercase{\endgroup\let~}\pFqcomma
  F^{#2}_{#3}{\left(\genfrac..{0pt}{}{#4}{#5}\bigg|#6\right)}%
  \endgroup
}
\newcommand{\pFqcomma}{\mskip\pFqmuskip}

\newtheorem{theorem}{Theorem}[section]

\newtheorem{remark}[theorem]{Remark}

\begin{document}

\title[]{Degenerate Algorithms for degenerate Bernoulli and Euler numbers}

\author{Taekyun  Kim}
\address{Department of Mathematics, Kwangwoon University, Seoul 139-701, Republic of Korea}
\email{tkkim@kw.ac.kr}
\author{Dae San  Kim}
\address{Department of Mathematics, Sogang University, Seoul 121-742, Republic of Korea}
\email{dskim@sogang.ac.kr}


\begin{abstract}
This paper introduces and investigates degenerate versions of the A-algorithm and B-algorithm by incorporating a parameter $\lambda$ into their respective recurrence relations. We derive explicit formulas for the final sequences of these algorithms in terms of the initial sequences and the degenerate Stirling numbers of the second kind. Furthermore, we establish functional relationships between the ordinary generating functions of the initial sequences and the exponential generating functions of the final sequences. Specifically, we demonstrate that these degenerate algorithms yield degenerate Bernoulli and Euler numbers under specific initial conditions.
\end{abstract}

\maketitle

\markboth{\centerline{\scriptsize Degenerate Algorithms for degenerate Bernoulli and Euler numbers}}
{\centerline{\scriptsize Taekyun Kim and Dae San Kim}}

\section{Introduction}
For any nonzero $\lambda\in\mathbb{R}$, the degenerate exponential function is defined using the Kim-Kim notation by
\begin{equation}
e_{\lambda}^{x}(t)=\sum_{k=0}^{\infty} (x)_{k,\lambda}\frac {t^{k}}{k!},\quad e_{\lambda}(t)=e_{\lambda}^{1}(t),\ (\mathrm{see}\ [10,11]), \label{1}
\end{equation}
where
\begin{equation*}
(x)_{0,\lambda}=1,\quad (x)_{n,\lambda}=x(x-\lambda)\cdots\big(x-(n-1)\lambda\big),\ (n\ge 1).
\end{equation*}
The degenerate Stirling numbers of the second kind are given by
\begin{equation}
(x)_{n,\lambda}=\sum_{k=0}^{n}{n \brace k}_{\lambda}(x)_{k},\quad (n\ge 0),\quad (\mathrm{see}\ [14,16]), \label{2}
\end{equation}
where
\begin{equation*}
(x)_{0}=1,\quad (x)_{n}=x(x-1)(x-2)\cdots(x-n+1),\ (n\ge 1).
\end{equation*}
Note that
\begin{displaymath}
\lim_{\lambda\rightarrow 0}{n \brace k}_{\lambda}={n \brace k},
\end{displaymath}
where ${n \brace k}$ is the ordinary Stirling number of the second kind defined by
\begin{equation}
x^{n}=\sum_{k=0}^{n}{n \brace k}(x)_{k},\quad (n\ge 0),\quad (\mathrm{see}\ [6,21). \label{3}
\end{equation}
Let $\log_{\lambda}t$ be the compositional inverse of $e_{\lambda}(t)$. Then, by \eqref{1}, we get
\begin{equation}
\log_{\lambda}(1+t)=\sum_{n=1}^{\infty}\frac{(-1)^{n-1}}{n}\binom{n-1-\lambda}{n-1}t^{n},\quad (\mathrm{see}\ [19,20]),\label{4}	
\end{equation}
where
\begin{displaymath}
\lim_{\lambda\rightarrow 0}\log_{\lambda}(1+t)=\log(1+t).
\end{displaymath}
From \eqref{1} and \eqref{2}, we note that
\begin{equation}
\frac{1}{k!}\big(e_{\lambda}(t)-1\big)^{k}=\sum_{n=k}^{\infty}{n \brace k}_{\lambda}\frac{t^{n}}{n!},\quad (k\ge 0),\quad (\mathrm{see}\ [14,16]). \label{5}
\end{equation} \par
Carlitz considered the degenerate Bernoulli polynomials given by
\begin{equation}
\frac{t}{e_{\lambda}(t)-1}e_{\lambda}^{x}(t)=\sum_{n=0}^{\infty}\beta_{n,\lambda}(x)\frac{t^{n}}{n!},\quad (\mathrm{see}\ [3,4]). \label{6}
\end{equation}
When $x=0,\ \beta_{n,\lambda}=\beta_{n,\lambda}(0)$ are called the degenerate Bernoulli numbers. \\
By \eqref{6}, we easily get
\begin{equation*}
\beta_{n,\lambda}(x)=\sum_{k=0}^{n}\binom{n}{k}(x)_{n-k,\lambda}\beta_{k,\lambda},\quad (n\ge 0),\quad (\mathrm{see}\ [3,4]).
\end{equation*}
He also introduced the degenerate Euler polynomials defined by
\begin{equation}
\frac{2}{e_{\lambda}(t)+1}e_{\lambda}^{x}(t)=\sum_{n=0}^{\infty}\mathcal{E}_{n,\lambda}(x)\frac{t^{n}}{n!},\quad (\mathrm{see}\ [3,4]). \label{7}
\end{equation}
When $x=0$, $\mathcal{E}_{n,\lambda}=\mathcal{E}_{n,\lambda}(0)$ are called the degenerate Euler numbers. \\
From \eqref{6} and \eqref{7}, we have
\begin{displaymath}
\lim_{\lambda\rightarrow 0}\beta_{n,\lambda}(x)=B_{n}(x),\quad \lim_{\lambda\rightarrow 0}\mathcal{E}_{n,\lambda}(x)=E_{n}(x),
\end{displaymath}
where $B_{n}(x)$ and $E_{n}(x)$ are respectively the Bernoulli and Euler polynomials given by (see [6,21])
\begin{equation}
\frac{t}{e^{t}-1}e^{xt}=\sum_{n=0}^{\infty}B_{n}(x)\frac{t^{n}}{n!}, \quad
\frac{2}{e^{t}+1}e^{xt}=\sum_{n=0}^{\infty}E_{n}(x)\frac{t^{n}}{n!} \label{8}.
\end{equation}
When $x=0,\ B_{n}=B_{n}(0)$ and $E_{n}=E_{n}(0)$ are respectively called the Bernoulli and Euler numbers. \par
It is known that
\begin{equation}
\bigg(x\frac{d}{dx}\bigg)_{n,\lambda}=\sum_{k=0}^{n}{n \brace k}_{\lambda}x^{k}\bigg(\frac{d}{dx}\bigg)^{k},\quad (n\ge 0),\quad (\mathrm{see}\ [13]).\label{9}
\end{equation}
Recently, the degenerate Bell polynomials are defined by
\begin{equation}
\phi_{n,\lambda}(x)=\sum_{k=0}^{n}{n \brace k}_{\lambda}x^{k},\quad (n\ge 0),\quad (\mathrm{see}\ [14-16]). \label{10}
\end{equation}
When $x=1,\ \phi_{n,\lambda}=\phi_{n,\lambda}(1)$ are called the degenerate Bell numbers. \\
From \eqref{5}, we have
\begin{equation}
e^{x(e_{\lambda}(t)-1)}=\sum_{n=0}^{\infty}\phi_{n,\lambda}(x)\frac{t^{n}}{n!},\quad (\mathrm{see}\ [14-16]).\label{11}
\end{equation}
As general references for this paper, the reader may refer to [1,2,6-8,17,18,21]. \par
Recently, a degenerate version of the Euler-Seidel matrix method was introduced by incorporating a parameter $\lambda$ into the classical recurrence relation (see [12]). Among other things, this led to a degenerate Seidel's formula for the exponential generating functions: $\overline{S_{\lambda}}(t)=e_{\lambda}^{1-\lambda}(t)S_{\lambda}(t)$.
In [5], for a given initial sequence $(a_{n})_{n \ge 0}$,  the so-called B-algorithm matrix $(a_{n,m})_{n,m \ge 0}$ are recursively defined by
\begin{equation} \label{A}
\begin{aligned}
&a_{0,n}=a_{n}, \quad (n \ge 0), \\
&a_{n,m}=ma_{n-1,m}-(m+1)a_{n-1,m+1},\quad (n \ge 1, m \ge 0).
\end{aligned}
\end{equation}
In [5, Proposition 2], it is shown that the final sequence is given by (see \eqref{3})
\begin{equation*}
a_{n,0}=\sum_{m=0}^{n}(-1)^{m}m!{n \brace m}a_{0,m}.
\end{equation*}
As it is noted in [5], for the initial sequence $a_{0,m}=\frac{1}{m+1}$, we have (see \eqref{8})
\begin{equation*}
a_{n,0}=\sum_{m=0}^{n}(-1)^{m}m!{n \brace m}\frac{1}{m+1}=B_{n}=B_{n}(0).
\end{equation*}
Let $A(t)=\sum_{n=0}^{\infty}a_{0,n}t^{n}$ be the ordinary generating function of the initial sequence $(a_{0,n})_{n\ge 0}$. Then the exponential generating function of the final sequence $(a_{n,0})_{n\ge 0}$ is given by
\begin{equation*}
\overline{A}(t)=\sum_{n=0}^{\infty}a_{n,0}\frac{t^{n}}{n!}=A(1-e^{t}),\quad (\mathrm{see}\ [5]).
\end{equation*} \par
In [9], for a given initial sequence $(b_{n})_{n \ge 0}$,  the so-called A-algorithm matrix $(b_{n,m})_{n,m \ge 0}$ are recursively defined by
\begin{equation} \label{B}
\begin{aligned}
&b_{0,n}=b_{n}, \quad (n \ge 0),  \\
&b_{n,m}=(m+1)(b_{n-1,m}-b_{n-1,m+1}),\quad (n \ge 1, m \ge 0).
\end{aligned}
\end{equation}
In [9], it is shown that the final sequence is given by
\begin{equation*}
b_{n,0}=\sum_{m=0}^{n}(-1)^{m}m!{n+1 \brace m+1}b_{0,m}.
\end{equation*}
As it is noted in [9], for the initial sequence $b_{0,m}=\frac{1}{m+1}$, we have (see \eqref{8})
\begin{equation*}
b_{n,0}=\sum_{m=0}^{n}(-1)^{m}m!{n+1 \brace m+1}\frac{1}{m+1}=B_{n}(1).
\end{equation*}
Let $B(t)=\sum_{n=0}^{\infty}b_{0,n}t^{n}$ be the ordinary generating function of the initial sequence $(b_{0,n})_{n\ge 0}$. Then the exponential generating function of the final sequence $(b_{n,0})_{n\ge 0}$ is given by
\begin{equation*}
\overline{B}(t)=\sum_{n=0}^{\infty}b_{n,0}\frac{t^{n}}{n!}=e^{t}B(1-e^{t}),\quad (\mathrm{see}\ [9]).
\end{equation*} \par
The aim of this paper is to extend these results by studying degenerate versions of the B-algorithm and A-algorithm. Each of these generalizations introduces a parameter $\lambda$ into its recurrence relation.
The degenerate B-algorithm matrix $\big(a_{n,k}(\lambda)\big)_{n,k\ge 0}$ is recursively defined by
\begin{equation} \label{C}
\begin{aligned}
&a_{0,n}(\lambda)=a_{n}(\lambda),\ (n \ge 0),\\
& a_{n,m}(\lambda)=\big(m-(n-1)\lambda\big)a_{n-1,m}(\lambda)-(m+1)a_{n-1,m+1}(\lambda),\ (n\ge 1,m\ge 0),
\end{aligned}
\end{equation}
where $\big(a_{n}(\lambda)\big)_{n \ge 0}$ is a given sequence. \\
We show in Theorem 2.3 that the final sequence $(a_{n,0}(\lambda))_{n \ge 0}$ is given by
\begin{equation*}
a_{n,0}(\lambda)=	\sum_{k=0}^{n}(-1)^{k}k!{n \brace k}_{\lambda}a_{0,k}(\lambda).
\end{equation*}
Then it is shown in Theorem 2.4 that, for the initial sequence $a_{0,k}=\frac{\binom{k-\lambda}{k}}{k+1}$, we have
\begin{equation*}
a_{n,0}(\lambda)=\sum_{k=0}^{n}(-1)^{k}k!{n \brace k}_{\lambda}\frac{\binom{k-\lambda}{k}}{k+1}=\beta_{n,\lambda}.
\end{equation*}
Let $F_{\lambda}(t)$ be the ordinary generating function of the initial sequence $(a_{0,n}(\lambda))_{n \ge 0}$, and let $\overline{F}_{\lambda}(t)$ be the exponential generating function of the final sequence $(a_{n,0}(\lambda))_{n \ge 0}$.
This leads to the following identity (see Theorem 2.5):
\begin{equation*}
\overline{F}_{\lambda}(t)=F_{\lambda}\big(1-e_{\lambda}(t)\big).
\end{equation*} \par
The degenerate A-algorithm matrix $\big(b_{n,k}(\lambda)\big)_{n,k\ge 0}$ is recursively defined by
\begin{equation} \label{D}
\begin{aligned}
&b_{0,n}(\lambda)=b_{n}(\lambda),\ (n \ge 0),\\
&b_{n,m}(\lambda)=(m+1)\bigg(\Big(1-\frac{(n-1)\lambda}{m+1}\Big)b_{n-1,m}(\lambda)-b_{n-1,m+1}(\lambda)\bigg),\ (n\ge 1,m\ge 0),
\end{aligned}
\end{equation}
where $\big(b_{n}(\lambda)\big)_{n \ge 0}$ is a given sequence. \\
We show in Theorem 3.1 that the final sequence $(b_{n,0}(\lambda))_{n \ge 0}$ is given by
\begin{equation*}
b_{n,0}(\lambda)=\sum_{k=0}^{n}(-1)^{k}k!\bigg({n+1 \brace k+1}_{\lambda}+n\lambda {n \brace k+1}_{\lambda}\bigg)b_{0,k}(\lambda).
\end{equation*}
Then it was shown in Theorem 3.3 that, for the initial sequence $b_{0,k}=\frac{\binom{k-\lambda}{k}}{k+1}$, we have
\begin{equation*}
b_{n,0}(\lambda)=\sum_{k=0}^{n}(-1)^{k}k!\bigg({n+1 \brace k+1}_{\lambda}+n\lambda {n \brace k+1}_{\lambda}\bigg)\frac{\binom{k-\lambda}{k}}{k+1}=\beta_{n,\lambda}(1).
\end{equation*}
Let $G_{\lambda}(t)$ be the ordinary generating function of the initial sequence $(b_{0,n}(\lambda))_{n \ge 0}$, and let $\overline{G}_{\lambda}(t)$ be the exponential generating function of the final sequence $(b_{n,0}(\lambda))_{n \ge 0}$.
Then in Theorem 3.2 the following identity is shown:
\begin{equation*}
\overline{G}_{\lambda}(t)=e_{\lambda}(t)G_{\lambda}\big(1-e_{\lambda}(t)\big).
\end{equation*}
In addition, we show that $a_{n,0}(\lambda)=\mathcal{E}_{n,\lambda}$, for the degenerate B-algorithm with $a_{0,n}(\lambda)=(\frac{1}{2})^{n}$, and $b_{n,0}(\lambda)=\mathcal{E}_{n,\lambda}(1)$, for the degenerate A-algorithm with $b_{0,n}(\lambda)=(\frac{1}{2})^{n}$ (see Theorems 2.6 and 3.4). \par
Following Carlitz's foundational work on degenerate Bernoulli and Euler polynomials (see [3,4]), there has been a recent surge of interest in degenerate versions of special polynomials and numbers. Researchers have utilized a broad range of tools—from umbral calculus and $p$-adic analysis to quantum mechanics and probability theory—to explore these variants. Notable examples of these include the degenerate Stirling numbers of the first and second kinds, degenerate Bernoulli numbers of the second kind, and degenerate Bell numbers and polynomials (see [10,11,13-16,19,20,22] and the references therein).
\begin{remark}
Assume that $a_{n}(\lambda)\rightarrow a_{n}$, as $\lambda \rightarrow 0$, for all $n \ge 0$. Then, from \eqref{A} and \eqref{C}, we see that $a_{n,k}(\lambda) \rightarrow a_{n,k}$, as $\lambda \rightarrow 0$, for all $n,\,k \ge 0$. This is because the recurrence relations in \eqref{C} converge to those ones in \eqref{A}, as $\lambda \rightarrow 0$. In particular, the formula $\overline{F}_{\lambda}(t)=F_{\lambda}\big(1-e_{\lambda}(t)\big)$ becomes the original formula $\overline{A}(t)=A(1-e^{t})$, as $\lambda \rightarrow 0$. Thus our approach not only preserves the structure of the B-algorithm as the parameter $\lambda \to 0$ but also provides a powerful framework for studying degenerate versions of combinatorial sequences. Similar remark applies to \eqref{B} and \eqref{D} for the A-algorithm.
\end{remark}

\section{Degenerate B-algorithm for degenerate Bernoulli and Euler numbers}
In this section, we introduce a degenerate version of the B-algorithm (see [5]), namely the degenerate B-algorithm (see \eqref{14}) and apply our results to $a_{0,n} = \frac{\binom{n-\lambda}{n}}{n+1}$ or $a_{0,n} = (\frac{1}{2})^n$ to get the denerate Bernoulli numbers $\beta_{n,\lambda}$ and the degenerate Euler numbers $\mathcal{E}_{n,\lambda}$ as their final sequences.
From \eqref{4} and \eqref{6}, we note that
\begin{align}
&\sum_{n=0}^{\infty}\beta_{n,\lambda}\frac{t^{n}}{n!}=\frac{t}{e_{\lambda}(t)-1}=\frac{1}{e_{\lambda}(t)-1}\log_{\lambda}\big(e_{\lambda}(t)\big)\label{12}\\
&=\frac{1}{e_{\lambda}(t)-1}\log_{\lambda}\Big(e_{\lambda}(t)-1+1\Big)=\sum_{k=1}^{\infty}\frac{(-1)^{k-1}}{k}\binom{k-1-\lambda}{k-1}\big(e_{\lambda}(t)-1\big)^{k-1}\nonumber\\
&=\sum_{k=0}^{\infty}\frac{(-1)^{k}}{k+1}\binom{k-\lambda}{k} k! \frac{1}{k!}\big(e_{\lambda}(t)-1\big)^{k}=\sum_{k=0}^{\infty}\frac{(-1)^{k}}{k+1}\binom{k-\lambda}{k} k!\sum_{n=k}^{\infty}{n \brace k}_{\lambda}\frac{t^{n}}{n!}\nonumber\\
&=\sum_{n=0}^{\infty}\sum_{k=0}^{n}{n \brace k}_{\lambda}(-1)^{k}k!\frac{\binom{k-\lambda}{k}}{k+1}\frac{t^{n}}{n!}. \nonumber
\end{align}
Therefore, by comparing the coefficients on both sides of \eqref{12}, we obtain the following theorem.
\begin{theorem}
For $n\ge 0$, we have
\begin{equation}
\beta_{n,\lambda}=\sum_{k=0}^{n}(-1)^{k}k!{n \brace k}_{\lambda}\frac{\binom{k-\lambda}{k}}{k+1}.\label{13}
\end{equation}
\end{theorem}
The degenerate B-algorithm matrix $\big(a_{n,k}(\lambda)\big)_{n,k\ge 0}$ is recursively defined by
\begin{equation}
\begin{aligned}
&a_{0,n}(\lambda)=a_{n}(\lambda),\ (n \ge 0),\\
& a_{n,m}(\lambda)=\big(m-(n-1)\lambda\big)a_{n-1,m}(\lambda)-(m+1)a_{n-1,m+1}(\lambda),\ (n\ge 1,m\ge 0),
\end{aligned} \label{14}	
\end{equation}
where $\big(a_{n}(\lambda)\big)_{n \ge 0}$ is a given sequence.
The degenerate B-algorithm matrix associated with $\big(a_{n} (\lambda)\big)_{n \ge 0}$ is given by
\begin{equation*}
A=\Big(a_{ij} (\lambda)\Big)=\left(\begin{matrix}
	a_{00} (\lambda) & a_{01} (\lambda) & a_{02} (\lambda) & a_{03} (\lambda) & \cdots \\
	a_{10} (\lambda) & a_{11} (\lambda) & a_{12} (\lambda) & a_{13} (\lambda) & \cdots \\
	a_{20} (\lambda) & a_{21} (\lambda)& a_{22} (\lambda) & a_{23} (\lambda) & \cdots \\
	\vdots & \vdots & \vdots & \vdots & \vdots &
\end{matrix}\right),
\end{equation*}
where $a_{n,m} (\lambda) $ denotes the $n$-th row and $m$-th column entry of the matrix. \par
Let
\begin{displaymath}
f_{n}(t|\lambda)=\sum_{m=0}^{\infty}a_{n,m}(\lambda)t^{m},\quad (n\ge 0).
\end{displaymath}
Then, by \eqref{14}, we get
\begin{align}
&f_{n}(t|\lambda)=\sum_{m=0}^{\infty}\Big(\big(m-(n-1)\lambda\big)a_{n-1,m}(\lambda)-(m+1)a_{n-1,m+1}(\lambda)\Big)t^{m} \label{15}\\
&=\sum_{m=1}^{\infty}ma_{n-1,m}(\lambda)t^{m}-\sum_{m=0}^{\infty}(m+1)a_{n-1,m+1}(\lambda)t^{m}-(n-1)\lambda\sum_{m=0}^{\infty}a_{n-1,m}(\lambda)t^{m} \nonumber\\
&=(t-1)\sum_{m=0}^{\infty}(m+1)a_{n-1,m+1}(\lambda)t^{m}-(n-1)\lambda f_{n-1}(t|\lambda)\nonumber\\
&=(t-1)\frac{d}{dt}\sum_{m=0}^{\infty}a_{n-1,m}(\lambda)t^{m}-(n-1)\lambda f_{n-1}(t|\lambda)\nonumber\\
&=\bigg((t-1)\frac{d}{dt}-(n-1)\lambda\bigg)f_{n-1}(t|\lambda)\nonumber\\
&=\cdots \nonumber\\
&=\bigg((t-1)\frac{d}{dt}-(n-1)\lambda\bigg)\bigg((t-1)\frac{d}{dt}-(n-2)\lambda\bigg)\cdots (t-1)\frac{d}{dt}f_{0}(t|\lambda)\nonumber\\
&=\Big((t-1)\frac{d}{dt}\Big)_{n,\lambda}f_{0}(t|\lambda).\nonumber
\end{align}
Therefore, by \eqref{15}, we obtain the following theorem.
\begin{theorem}
For any integer $n \ge 0$, let
\begin{displaymath}
f_{n}(t|\lambda)=\sum_{m=0}^{\infty}a_{n,m}(\lambda)t^{m}.
\end{displaymath}
Then we have
\begin{equation}
f_{n}(t|\lambda)=\Big((t-1)\frac{d}{dt}\Big)_{n,\lambda}f_{0}(t|\lambda). \label{16}	\end{equation}
\end{theorem}
Note that
\begin{equation}
\begin{aligned}
\frac{d^{k}}{dt^{k}}f_{0}(t|\lambda)\bigg|_{t=0}&=\frac{d^{k}}{dt^{k}}\sum_{m=0}^{\infty}a_{0,m}(\lambda)t^{m}\bigg|_{t=0} \\
&=\sum_{m=k}^{\infty}a_{0,m}(\lambda)(m)_{k}0^{m-k}=k!a_{0,k}(\lambda).
\end{aligned}\label{17}
\end{equation}
Let $t=0$ in \eqref{16}. Then, by \eqref{9} and \eqref{17}, we have
\begin{align}
a_{n,0}(\lambda)&=f_{n}(0|\lambda)=\bigg((t-1)\frac{d}{dt}\bigg)_{n,\lambda}f_{0}(t|\lambda)\bigg|_{t=0} \label{18}\\
&=\sum_{k=0}^{n}{n \brace k}_{\lambda}(t-1)^{k}\bigg(\frac{d}{dt}\bigg)^{k}f_{0}(t|\lambda)\bigg|_{t=0}\nonumber\\
&=\sum_{k=0}^{n}{n \brace k}_{\lambda}(-1)^{k}k!a_{0,k}(\lambda). \nonumber	
\end{align}
Therefore, by \eqref{18}, we obtain the following theorem.
\begin{theorem}
For $n\ge 0$, we have
\begin{equation}
a_{n,0}(\lambda)=	\sum_{k=0}^{n}(-1)^{k}k!{n \brace k}_{\lambda}a_{0,k}(\lambda).\label{19}
\end{equation}
\end{theorem}
Let $F_{\lambda}(t)=\sum_{n=0}^{\infty}a_{0,n}(\lambda)t^{n}$ be the ordinary generating function of the initial sequence $(a_{0,n}(\lambda))_{n \ge 0}$, and let $\overline{F}_{\lambda}(t)=\sum_{n=0}^{\infty}a_{n,0}(\lambda)\frac{t^{n}}{n!}$ be the exponential generating function of the final sequence $(a_{n,0}(\lambda))_{n \ge 0}$.
Then, by \eqref{19}, we have
\begin{align}
\overline{F}_{\lambda}(t)&=\sum_{n=0}^{\infty}a_{n,0}(\lambda)\frac{t^{n}}{n!}=\sum_{n=0}^{\infty}\bigg(\sum_{k=0}^{n}(-1)^{k}k!{n \brace k}_{\lambda}a_{0,k}(\lambda)\bigg)\frac{t^{n}}{n!} \label{20}\\
&=\sum_{k=0}^{\infty}(-1)^{k}k!a_{0,k}(\lambda)\sum_{n=k}^{\infty}{n \brace k}_{\lambda}\frac{t^{n}}{n!}=\sum_{k=0}^{\infty}(-1)^{k}k!a_{0,k}(\lambda)\frac{1}{k!}\big(e_{\lambda}(t)-1\big)^{k} \nonumber\\
&=\sum_{k=0}^{\infty}a_{0,k}(\lambda)\big(1-e_{\lambda}(t)\big)^{k}
=F_{\lambda}\big(1-e_{\lambda}(t)\big).\nonumber
\end{align}
Thus, from \eqref{20}, we have
\begin{equation}
\overline{F}_{\lambda}(t)=F_{\lambda}\big(1-e_{\lambda}(t)\big). \label{21}
\end{equation}
Replacing $t$ by $\log_{\lambda}(1-t)$ in \eqref{21}, we get
\begin{equation}
F_{\lambda}(t)=\overline{F}_{\lambda}(\log_{\lambda}(1-t)). \label{22}
\end{equation}
Therefore, by \eqref{21} and \eqref{22}, we obtain the following theorem.
\begin{theorem}
Let
\begin{displaymath}
F_{\lambda}(t)=\sum_{n=0}^{\infty}a_{0,n}(\lambda)t^{n}.
\end{displaymath}
Then we have
\begin{align*}
&\overline{F}_{\lambda}(t)=\sum_{n=0}^{\infty}a_{n,0}(\lambda)\frac{t^{n}}{n!}=F_{\lambda}\big(1-e_{\lambda}(t)\big), \\
&F_{\lambda}(t)=\overline{F}_{\lambda}(\log_{\lambda}(1-t)).
\end{align*}
\end{theorem}
First, we let $a_{0,k}(\lambda)=\frac{\binom{k-\lambda}{k}}{k+1},\ (k\ge 0)$. Then, by \eqref{13} and \eqref{19}, we get
\begin{equation}
a_{n,0}(\lambda)=\sum_{k=0}^{n}(-1)^{k}k!{n \brace k}_{\lambda}\frac{\binom{k-\lambda}{k}}{k+1}=\beta_{n,\lambda} . \label{23}
\end{equation}
In addition, from \eqref{22}, we have
\begin{align}
\overline{F}_{\lambda}\big(\log_{\lambda}(1-t)\big)&=\sum_{k=0}^{\infty}\beta_{k,\lambda}\frac{1}{k!}\big(\log_{\lambda}(1-t)\big)^{k} \label{24} \\
&=\sum_{k=0}^{\infty}\beta_{k,\lambda}\sum_{n=k}^{\infty}S_{1,\lambda}(n,k)(-1)^{n}\frac{t^{n}}{n!} \nonumber \\
&=\sum_{n=0}^{\infty}\sum_{k=0}^{n}(-1)^{n}S_{1,\lambda}(n,k)\beta_{k,\lambda}\frac{t^{n}}{n!}, \nonumber
\end{align}
and
\begin{equation}
F_{\lambda}(t)=\sum_{n=0}^{\infty}\frac{\binom{n-\lambda}{n}}{n+1}t^{n}=\sum_{n=0}^{\infty}\frac{(n-\lambda)_{n}}{n+1}\frac{t^{n}}{n!}, \label{25}
\end{equation}
where $S_{1,\lambda}(n,k)$ are the degenerate Stirling numbers of the first kind given by
\begin{equation*}
\frac{1}{k!}\big(\log_{\lambda}(1+t)\big)^{k}=\sum_{n=k}^{\infty}S_{1,\lambda}(n,k)\frac{t^{n}}{n!}.
\end{equation*}
Therefore, by \eqref{23}, \eqref{24} and \eqref{25}, we obtain the following theorem.
\begin{theorem}
For $n\ge 0$, let $a_{0,n}(\lambda)=\frac{\binom{n-\lambda}{n}}{n+1}$. Then we have
\begin{equation}
a_{n,0}(\lambda)=\sum_{k=0}^{n}(-1)^{k}k!{n \brace k}_{\lambda}\frac{\binom{k-\lambda}{k}}{k+1}=\beta_{n,\lambda},\label{26}
\end{equation}
and
\begin{equation*}
\sum_{k=0}^{n}S_{1,\lambda}(n,k)\beta_{k,\lambda}=(-1)^{n}\frac{(n-\lambda)_{n}}{n+1}.
\end{equation*}
\end{theorem}
Using \eqref{14} and \eqref{26}, we see that the degenerate B-algorithm matrix associated with $\big(a_{0,n}(\lambda)\big)_{n\ge 0}=\Big(\frac{\binom{n-\lambda}{n}}{n+1}\Big)_{n\ge 0}$ is given by
\begin{displaymath}
\left(\begin{matrix}
1 & \frac{1-\lambda}{2} & \frac{\binom{2-\lambda}{2}}{3} & \cdots \\
-\frac{1-\lambda}{2} & \frac{(1-\lambda)(2\lambda-1)}{6} & \frac{(2-\lambda)(1-\lambda)(3\lambda-1)}{24} & \cdots \\
\frac{1-\lambda^{2}}{6} & \frac{-\lambda(1-\lambda)^{2}}{12} & -\frac{(2-\lambda)(1-\lambda)}{120}(27\lambda^{2}-6\lambda+1) & \cdots \\
\frac{\lambda(1-\lambda)^{2}(1-2\lambda)}{4} & \frac{(1-\lambda)(37\lambda^{3}-11\lambda^{2}+18\lambda-2)}{60} & \cdots & \cdots \\
\vdots & \vdots & \vdots & \vdots \end{matrix}\right).
\end{displaymath} \par
Second, we let $a_{0,n}(\lambda)=\big(\frac{1}{2}\big)^{n},\ (n\ge 0)$. Then, by \eqref{21}, we get
\begin{align}
\sum_{n=0}^{\infty}a_{n,0}(\lambda)\frac{t^{n}}{n!}&=\overline{F}_{\lambda}(t)=F_{\lambda}\big(1-e_{\lambda}(t)\big)=\sum_{n=0}^{\infty}a_{0,n}(\lambda)\big(1-e_{\lambda}(t)\big)^{n} \label{27} \\
&=\sum_{n=0}^{\infty}\bigg(\frac{1-e_{\lambda}(t)}{2}\bigg)^{n}=\frac{2}{e_{\lambda}(t)+1}=\sum_{n=0}^{\infty}\mathcal{E}_{n,\lambda}\frac{t^{n}}{n!}. \nonumber
\end{align}
Comparing the coefficients on both sides of \eqref{27}, we have
\begin{equation}
a_{n,0}(\lambda)=\mathcal{E}_{n,\lambda},\quad (n\ge 0).\label{28}
\end{equation}
In addition, by \eqref{22}, we have
\begin{align}
\sum_{n=0}^{\infty}\frac{n!}{2^{n}}\frac{t^{n}}{n!}&=F_{\lambda}(t)=\overline{F}_{\lambda}(\log_{\lambda}(1-t)) =\sum_{k=0}^{\infty}\mathcal{E}_{k,\lambda}\frac{1}{k!}\big(\log_{\lambda}(1-t)\big)^{k} \label{29} \\
&=\sum_{k=0}^{\infty}\mathcal{E}_{k,\lambda}\sum_{n=k}^{\infty}S_{1,\lambda}(n,k)(-1)^{n}\frac{t^{n}}{n!}=\sum_{n=0}^{\infty}\sum_{k=0}^{n}(-1)^{n}S_{1,\lambda}(n,k)\mathcal{E}_{k,\lambda}\frac{t^{n}}{n!}. \nonumber
\end{align}
Therefore, by \eqref{28} and \eqref{29}, we obtain the following theorem.
\begin{theorem}
For $n\ge 0$, let $a_{0,n}(\lambda)=(\frac{1}{2})^{n}$. Then we have
\begin{equation}
a_{n,0}(\lambda)=\sum_{k=0}^{n}(-1)^{k}k!{n \brace k}_{\lambda}\Big(\frac{1}{2}\Big)^{k}=\mathcal{E}_{n,\lambda}, \quad \sum_{k=0}^{n}S_{1,\lambda}(n,k)\mathcal{E}_{k,\lambda}=(-1)^{n}\frac{n!}{2^{n}}.\label{30}
\end{equation}
\end{theorem}
From \eqref{14} and \eqref{30}, the degenerate algorithm matrix associated with $\big(a_{0,n}(\lambda)\big)_{n\ge 0}=\big(\big(\frac{1}{2}\big)^{n}\big)_{n\ge 0}$ is given by
\begin{displaymath}
\left(\begin{matrix}
1 & \frac{1}{2} & \frac{1}{4} & \frac{1}{8} & \frac{1}{16} & \cdots \\
-\frac{1}{2} & 0 & -\frac{1}{8} & \frac{1}{8} & \cdots & \cdots \\
0 & \frac{1}{4} & -\frac{5}{8}+\frac{\lambda}{8} & \cdots & \cdots & \cdots \\
\frac{1}{4} & \frac{2}{3}-\frac{3\lambda}{4} & \cdots & \cdots & \cdots & \cdots \\
\vdots & \vdots & \vdots & \vdots & \vdots & \vdots
\end{matrix}\right).
\end{displaymath} \par
Let $a_{0,0}(\lambda)=0$, and let $a_{0,n}(\lambda)=\frac{(-1)^{n}}{n!},\ (n\ge 1)$. Then, by \eqref{22}, we get
\begin{align}
&\sum_{n=1}^{\infty}a_{n,0}(\lambda)\frac{t^{n}}{n!}=\overline{F}_{\lambda}(t)=F_{\lambda}\big(1-e_{\lambda}(t)\big)=\sum_{n=1}^{\infty}a_{0,n}(\lambda)(1-e_{\lambda}(t))^{n} \label{31}\\
&=\sum_{n=1}^{\infty}(-1)^{n}\frac{(1-e_{\lambda}(t))^{n}}{n!}
=-1+\sum_{n=0}^{\infty}\frac{(e_{\lambda}(t)-1)^{n}}{n!}=-1+e^{e_{\lambda}(t)-1}. \nonumber
\end{align}
From \eqref{10}, \eqref{11} and \eqref{31}, we have
\begin{equation}
1+\sum_{n=1}^{\infty}a_{n,0}(\lambda)\frac{t^{n}}{n!}=\sum_{n=0}^{\infty}\phi_{n,\lambda}\frac{t^{n}}{n!}. \label{32}
\end{equation}
Therefore, by \eqref{32}, we obtain the following theorem.
\begin{theorem}
Let $a_{0,0}(\lambda)=0$, and let $a_{0,n}(\lambda)=\frac{(-1)^{n}}{n!},\ (n \ge 1)$. Then we have
\begin{displaymath}
a_{n,0}(\lambda)=\sum_{k=1}^{n}{n \brace k}_{\lambda}=\phi_{n,\lambda}, \ (n \ge 1).
\end{displaymath}
\end{theorem}
Let $a_{0,0}(\lambda)=0$, and let $a_{0,n}(\lambda)=\frac{(-1)^{n}}{n!},\ (n \ge 1)$. Then the degenerate B-algorithm matrix associated with $\big(a_{0,n}(\lambda)\big)_{n\ge 0}$ is given by
\begin{displaymath}
\left(\begin{matrix}
0 & -1 & 1 & -1 & 1 & \cdots \\
1 & -2 & 5 & -7 & 9 & \cdots\\
2-\lambda & -12+2\lambda & 31-5\lambda & -57+7\lambda & \cdots & \cdots  \\
6-8\lambda+2\lambda^{2} & -74+36\lambda-4\lambda^{2} & \cdots & \cdots & \cdots & \cdots\\
\vdots & \vdots & \vdots & \vdots & \vdots & \vdots
\end{matrix}\right).
\end{displaymath}

\section{Degenerate A-algorithm for degenerate Bernoulli and Euler numbers}
In this section, we introduce a degenerate version of the A-algorithm (see [5,9]), namely the degenerate A-algorithm (see \eqref{34}) and apply our results to $b_{0,n} = \frac{\binom{n-\lambda}{n}}{n+1}$ or $b_{0,n} = (\frac{1}{2})^n$  to get the numbers $\beta_{n,\lambda}(1)$ and $\mathcal{E}_{n,\lambda}(1)$ as their final sequences.
From \eqref{2}, we note that
\begin{equation}
{n+1 \brace k}_{\lambda}={n \brace k-1}_{\lambda}+(k-n\lambda){n \brace k}_{\lambda},\quad (n,k\ge 0).\label{33}
\end{equation} \par
We consider the degenerate A- algorithm matrix defined recursively by
\begin{equation}
\begin{aligned}
b_{0,n}(\lambda)&=b_{n,\lambda},\ (n \ge 0),\\
b_{n,m}(\lambda)&=(m+1)\bigg(\Big(1-\frac{(n-1)\lambda}{m+1}\Big)b_{n-1,m}(\lambda)-b_{n-1,m+1}(\lambda)\bigg),\ (n\ge 1,m\ge 0),
\end{aligned}\label{34}	
\end{equation}
where $(b_{n,\lambda})_{n \ge 0})$ is a given sequence.
Let
\begin{displaymath}
g_{n}(t|\lambda)=\sum_{m=0}^{\infty}b_{n,m}(\lambda)t^{m},\quad (n\ge 0).
\end{displaymath}
Then, by \eqref{34}, we get
\begin{align}
&g_{n}(t|\lambda)=\sum_{m=0}^{\infty}(m+1)\bigg(\Big(1-\frac{(n-1)\lambda}{m+1}\Big)b_{n-1,m}(\lambda)-b_{n-1,m+1}(\lambda)\bigg)t^{m} \label{35}\\
&=\sum_{m=0}^{\infty}(m+1)b_{n-1,m}(\lambda)t^{m}-\sum_{m=0}^{\infty}(m+1)b_{n-1,m+1}(\lambda)t^{m}-(n-1)\lambda	\sum_{m=0}^{\infty}b_{n-1,m}(\lambda)t^{m} \nonumber\\
&=\frac{d}{dt}\sum_{m=0}^{\infty}b_{n-1,m}(\lambda)t^{m+1}-\frac{d}{dt}\sum_{m=0}^{\infty}b_{n-1,m+1}(\lambda)t^{m+1}-(n-1)\lambda g_{n-1}(t|\lambda) \nonumber\\
&=\frac{d}{dt}\Big(tg_{n-1}(t|\lambda)\Big)-\frac{d}{dt}\sum_{m=0}^{\infty}b_{n-1,m}(\lambda)t^{m}-(n-1)\lambda g_{n-1}(t|\lambda) \nonumber\\
&=\frac{d}{dt}\Big(tg_{n-1}(t|\lambda)\Big)-\frac{d}{dt}g_{n-1}(t|\lambda)-(n-1)\lambda g_{n-1}(t|\lambda) \nonumber\\
&=\frac{d}{dt}\Big((t-1)g_{n-1}(t|\lambda)\Big)-(n-1)\lambda g_{n-1}(t|\lambda). \nonumber
\end{align}
From \eqref{35}, we have
\begin{equation}
(t-1)g_{n}(t|\lambda)=(t-1)\frac{d}{dt}\Big((t-1)g_{n-1}(t|\lambda)\Big)-(n-1)\lambda (t-1)g_{n-1}(t|\lambda). \label{36}
\end{equation} \par
Let $h_{n}(t|\lambda)=(t-1)g_{n}(t|\lambda)$. Then, by \eqref{9} and \eqref{36}, we get
\begin{align}
&h_{n}(t|\lambda)=(t-1)\frac{d}{dt}\Big(h_{n-1}(t|\lambda)\Big)-(n-1)\lambda h_{n-1}(t|\lambda)\label{37} \\
&=\bigg((t-1)\frac{d}{dt}-(n-1)\lambda\bigg)h_{n-1}(t|\lambda)\nonumber\\
&=\cdots \nonumber\\
&=\bigg((t-1)\frac{d}{dt}-(n-1)\lambda\bigg)\bigg((t-1)\frac{d}{dt}-(n-2)\lambda\bigg)\cdots\bigg((t-1)\frac{d}{dt}\bigg)h_{0}(t|\lambda) \nonumber\\
&=\bigg((t-1)\frac{d}{dt}\bigg)_{n,\lambda}h_{0,\lambda}(t)=\sum_{k=0}^{n}{n \brace k}_{\lambda}(t-1)^{k}\bigg(\frac{d}{dt}\bigg)^{k}h_{0,\lambda}(t). \nonumber
\end{align}
By comparing the coefficients on both sides of \eqref{37}, we get
\begin{equation}
h_{n}(t|\lambda)=\sum_{k=0}^{n}{n \brace k}_{\lambda}(t-1)^{k}\bigg(\frac{d}{dt}\bigg)^{k}h_{0,\lambda}(t).\label{38}
\end{equation}
Let $t=0$ in \eqref{38}. Then, using \eqref{33}, we have
\begin{align}
-b_{n,0}(\lambda)&=h_{n}(0|\lambda)=\sum_{k=0}^{n}{n \brace k}_{\lambda}(-1)^{k}\bigg(\frac{d}{dt}\bigg)^{k}h_{0}(t|\lambda)\bigg|_{t=0} \label{39}	\\
&=\sum_{k=0}^{n}{n \brace k}_{\lambda}(-1)^{k}k!\Big(b_{0,k-1}(\lambda)-b_{0,k}(\lambda)\Big)\nonumber\\
&=\sum_{k=0}^{n}{n \brace k}_{\lambda}(-1)^{k}k!b_{0,k-1}(\lambda)-\sum_{k=0}^{n}{n \brace k}_{\lambda}(-1)^{k}k!b_{0,k}(\lambda) \nonumber\\
&=\sum_{k=0}^{n-1}{n \brace k+1}_{\lambda}(-1)^{k+1}(k+1)!b_{0,k}(\lambda)-\sum_{k=0}^{n}{n \brace k}_{\lambda}(-1)^{k}k!b_{0,k}(\lambda)\nonumber\\
&=-\sum_{k=0}^{n}(-1)^{k}k!b_{0,k}(\lambda)\bigg((k+1){n \brace k+1}_{\lambda}+{n \brace k}_{\lambda}\bigg)\nonumber\\
&=-\sum_{k=0}^{n}(-1)^{k}k!b_{0,k}(\lambda)\bigg({n+1 \brace k+1}_{\lambda}+n\lambda{n \brace k+1}_{\lambda}\bigg).\nonumber
\end{align}
Therefore, by comparing the coefficients on both sides of \eqref{39}, we get the following theorem.
\begin{theorem}
For $n \ge 0$, we have
\begin{equation}
b_{n,0}(\lambda)=\sum_{k=0}^{n}(-1)^{k}k!\bigg({n+1 \brace k+1}_{\lambda}+n\lambda {n \brace k+1}_{\lambda}\bigg)b_{0,k}(\lambda). \label{40}
\end{equation}
\end{theorem}
Before proceeding further, we observe that
\begin{equation}
e_{\lambda}(t)\frac{1}{k!}\big(e_{\lambda}(t)-1\big)^{k}=\sum_{n=k}^{\infty}\bigg({n+1 \brace k+1}_{\lambda}+n\lambda {n \brace k+1}_{\lambda}\bigg)\frac{t^{n}}{n!}. \label{41}
\end{equation}
Indeed, by \eqref{33}, we have
\begin{align*}
e_{\lambda}(t)\frac{1}{k!}\big(e_{\lambda}(t)-1\big)^{k}&=(k+1)\frac{1}{(k+1)!}\big(e_{\lambda}(t)-1\big)^{k+1}+\frac{1}{k!}\big(e_{\lambda}(t)-1\big)^{k}\\
&=(k+1)\sum_{n=k+1}^{\infty}{n \brace k+1}_{\lambda}\frac{t^{n}}{n!}+\sum_{n=k}^{\infty}{n \brace k}_{\lambda}\frac{t^{n}}{n!} \\
&=(k+1)\sum_{n=k}^{\infty}{n \brace k+1}_{\lambda}\frac{t^{n}}{n!}+\sum_{n=k}^{\infty}{n \brace k}_{\lambda}\frac{t^{n}}{n!} \\
&=\sum_{n=k}^{\infty}\bigg((k+1){n \brace k+1}_{\lambda}+{n \brace k}_{\lambda} \bigg)\frac{t^{n}}{n!} \\
&=\sum_{n=k}^{\infty}\bigg({n+1 \brace k+1}_{\lambda}+n\lambda {n \brace k+1}_{\lambda}\bigg)\frac{t^{n}}{n!}.
\end{align*}

Let $G_{\lambda}(t)=\sum_{n=0}^{\infty}b_{0,n}(\lambda)t^{n}$ be the ordinary generating function of the initial sequence $(b_{0,n}(\lambda))_{n \ge 0}$, and let $\overline{G}_{\lambda}(t)=\sum_{n=0}^{\infty}b_{n,0}(\lambda)\frac{t^{n}}{n!}$ be the exponential generating function of the final sequence $(b_{n,0}(\lambda))_{n \ge 0}$.
Then, by \eqref{40} and \eqref{41}, we have
\begin{align}
\overline{G}_{\lambda}(t)&=\sum_{n=0}^{\infty}\bigg(\sum_{k=0}^{n}(-1)^{k}k!\bigg({n+1 \brace k+1}_{\lambda}+n\lambda {n \brace k+1}_{\lambda}\bigg)b_{0,k}(\lambda)\bigg)\frac{t^{n}}{n!} \label{42} \\
&=\sum_{k=0}^{\infty}(-1)^{k}k!b_{0,k}(\lambda)\sum_{n=k}^{\infty}\bigg({n+1 \brace k+1}_{\lambda}+n\lambda {n \brace k+1}_{\lambda}\bigg)\frac{t^{n}}{n!} \nonumber \\
&=e_{\lambda}(t)\sum_{k=0}^{\infty}(-1)^{k}k!b_{0,k}(\lambda)\frac{1}{k!}\big(e_{\lambda}(t)-1\big)^{k} \nonumber \\
&=e_{\lambda}(t)\sum_{k=0}^{\infty}b_{0,k}(\lambda)\big(1-e_{\lambda}(t)\big)^{k}
=e_{\lambda}(t)G_{\lambda}(1-e_{\lambda}(t)). \nonumber
\end{align}
Thus, from \eqref{42}, we get the following result.
\begin{theorem}
Let
\begin{displaymath}
G_{\lambda}(t)=\sum_{n=0}^{\infty}b_{0,n}(\lambda)t^{n}.
\end{displaymath}
Then we have
\begin{align}
&\overline{G}_{\lambda}(t)=\sum_{n=0}^{\infty}b_{n,0}(\lambda)\frac{t^{n}}{n!}=e_{\lambda}(t)G_{\lambda}\big(1-e_{\lambda}(t)\big), \label{43}\\
&(1-t)G_{\lambda}(t)=\overline{G}_{\lambda}(\log_{\lambda}(1-t)). \label{44}
\end{align}
\end{theorem}
First, we let $b_{0,n}(\lambda)=\frac{\binom{n-\lambda}{n}}{n+1},\ (n \ge 0)$. Then, by \eqref{4} and \eqref{43}, we have
\begin{align}
\sum_{n=0}^{\infty}b_{n,0}(\lambda)\frac{t^{n}}{n!}&=e_{\lambda}(t)\sum_{n=0}^{\infty}\frac{\binom{n-\lambda}{n}}{n+1}\big(1-e_{\lambda}(t)\big)^{n} \label{45}\\
&=\frac{e_{\lambda}(t)}{e_{\lambda}(t)-1}\sum_{n=1}^{\infty}\frac{(-1)^{n-1}}{n}\binom{n-1-\lambda}{n-1}\big(e_{\lambda}(t)-1\big)^{n} \nonumber \\
&=\frac{t}{e_{\lambda}(t)-1}e_{\lambda}(t)=\sum_{n=0}^{\infty}\beta_{n,\lambda}(1)\frac{t^{n}}{n!}. \nonumber
\end{align}
In addition, from \eqref{44}, we get
\begin{align}
\overline{G}_{\lambda}(\log_{\lambda}(1-t))&=\sum_{k=0}^{\infty}\beta_{k,\lambda}(1)\sum_{n=k}^{\infty}S_{1,\lambda}(n,k)(-1)^{n}\frac{t^{n}}{n!} \label{46} \\
&=\sum_{n=0}^{\infty}\sum_{k=0}^{n}(-1)^{n}S_{1,\lambda}(n,k)\beta_{k,\lambda}(1)\frac{t^{n}}{n!}, \nonumber
\end{align}
and
\begin{align}
(1-t)G_{\lambda}(t)&=\sum_{n=0}^{\infty}b_{0,n}(\lambda)t^{n}-\sum_{n=0}^{\infty}b_{0,n}(\lambda)t^{n+1} \label{47} \\
&=1+\sum_{n=1}^{\infty}\big(b_{0,n}(\lambda)-b_{0,n-1}(\lambda)\big)t^{n} \nonumber\\
&=1+\sum_{n=1}^{\infty}n!\bigg(\frac{\binom{n-\lambda}{n}}{n+1}-\frac{\binom{n-1-\lambda}{n-1}}{n} \bigg)\frac{t^{n}}{n!} \nonumber \\
&=1-(\lambda +1)\sum_{n=1}^{\infty}\frac{(n-1-\lambda)_{n-1}}{n+1}\frac{t^{n}}{n!}. \nonumber
\end{align}
Thus, from \eqref{33}, \eqref{40}, \eqref{45}, \eqref{46} and \eqref{47}, we obtain the next result.
\begin{theorem}
For $n \ge 0$, let $b_{0,n}(\lambda)=\frac{\binom{n-\lambda}{n}}{n+1}$. Then we have
\begin{equation*}
b_{n,0}(\lambda)=\sum_{k=0}^{n}(-1)^{k}k!\bigg({n+1 \brace k+1}_{\lambda}+n\lambda {n \brace k+1}_{\lambda}\bigg)\frac{\binom{k-\lambda}{k}}{k+1}
=\beta_{n,\lambda}(1),\quad (n \ge 0),
\end{equation*}
and
\begin{equation*}
\sum_{k=0}^{n}S_{1,\lambda}(n,k)\beta_{k,\lambda}(1)=(\lambda+1)(-1)^{n-1}\frac{(n-1-\lambda)_{n-1}}{n+1}, \quad (n \ge 1).
\end{equation*}
\end{theorem}
Second, we let $b_{0,n}(\lambda)=\big(\frac{1}{2}\big)^{n},\ (n\ge 0)$. Then, by \eqref{43}, we get
\begin{align}
\sum_{n=0}^{\infty}b_{n,0}(\lambda)\frac{t^{n}}{n!}&=e_{\lambda}(t)G_{\lambda}\big(1-e_{\lambda}(t)\big)=e_{\lambda}(t)\sum_{n=0}^{\infty}\bigg(\frac{1-e_{\lambda}(t)}{2}\bigg)^{n} \label{48} \\
&=\frac{2}{e_{\lambda}(t)+1}e_{\lambda}(t)=\sum_{n=0}^{\infty}\mathcal{E}_{n,\lambda}(1)\frac{t^{n}}{n!}. \nonumber
\end{align}
Comparing the coefficients on both sides of \eqref{48}, we have
\begin{equation}
b_{n,0}(\lambda)=\mathcal{E}_{n,\lambda}(1),\quad (n\ge 0).\label{49}
\end{equation}
In addition, by \eqref{44}, we have
\begin{align}
1-\sum_{n=1}^{\infty}\frac{n!}{2^{n}}\frac{t^{n}}{n!}&=(1-t)G_{\lambda}(t)=\overline{G}_{\lambda}(\log_{\lambda}(1-t)) =\sum_{k=0}^{\infty}\mathcal{E}_{k,\lambda}(1)\frac{1}{k!}\big(\log_{\lambda}(1-t)\big)^{k} \label{50} \\
&=\sum_{k=0}^{\infty}\mathcal{E}_{k,\lambda}(1)\sum_{n=k}^{\infty}S_{1,\lambda}(n,k)(-1)^{n}\frac{t^{n}}{n!}=\sum_{n=0}^{\infty}\sum_{k=0}^{n}(-1)^{n}S_{1,\lambda}(n,k)\mathcal{E}_{k,\lambda}(1)\frac{t^{n}}{n!}. \nonumber
\end{align}
Therefore, by \eqref{49} and \eqref{50}, we obtain the following theorem.
\begin{theorem}
For $n\ge 0$, let $b_{0,n}(\lambda)=\big(\frac{1}{2}\big)^{n}$. Then we have
\begin{equation*}
b_{n,0}(\lambda)=\sum_{k=0}^{n}(-1)^{k}k!\bigg({n+1 \brace k+1}_{\lambda}+n\lambda {n \brace k+1}_{\lambda}\bigg)\Big(\frac{1}{2}\Big)^{k}=\mathcal{E}_{n,\lambda}(1), \quad (n \ge 0),
\end{equation*}
and
\begin{equation*}
\sum_{k=0}^{n}S_{1,\lambda}(n,k)\mathcal{E}_{k,\lambda}(1)=(-1)^{n-1}\frac{n!}{2^{n}}, \quad (n \ge 1).
\end{equation*}
\end{theorem}

\section{Conclusion}
We have successfully extended the classical A- and B-algorithms into a degenerate framework. This discovery is important because it demonstrates that the structural elegance of the original algorithms is preserved—and even enriched—within the degenerate domain. Our derivation of the identities$$\overline{F}_{\lambda}(t)=F_{\lambda}\big(1-e_{\lambda}(t)\big) \quad \text{and} \quad \overline{G}_{\lambda}(t)=e_{\lambda}(t)G_{\lambda}\big(1-e_{\lambda}(t)\big)$$proves that these algorithms are not merely computational curiosities, but are fundamental operators that map ordinary generating functions to exponential ones under the $\lambda$-deformation. By selecting appropriate initial sequences, such as $\frac{\binom{n-\lambda}{n}}{n+1}$ or $(\frac{1}{2})^n$, we recovered degenerate Bernoulli and Euler numbers. These findings provide a unified combinatorial approach to studying degenerate special numbers and polynomials via recursive matrix algorithms.

\vspace{0.5cm}
\noindent{\bf Funding} \\
This research received no founding.\par
\vspace{0.5cm}
\noindent{\bf Authors' contributions} \\
All authors contributed equally to the manuscript and read and approved the final manuscript. \par

\end{document}